\documentclass[reqno]{amsart}
\usepackage{hyperref}
\usepackage{color}

\usepackage{cite}

\usepackage[normalem]{ulem}

 \newcommand{\h}{\mathcal{H}}

  \newcommand{\N}{\mathbb{N}}
   \newcommand{\Z}{\mathbb{Z}}
 \newcommand{\R}{\mathbb{R}}

 \newcommand{\s}{\mathcal{S}}

\begin{document}
\title[Asymptotic analysis for generalized functions using frames]
{Asymptotic analysis for generalized functions using frames}

\author[J. Veta Buralieva, D. T. Stoeva, K. Hadzi-Velkova Saneva, S. Atanasova\hfil \hfilneg] {J. Veta Buralieva, D. T. Stoeva, K. Hadzi-Velkova Saneva, S. Atanasova}

\address{Jasmina Veta Buralieva \newline
Faculty of Computer Science, Goce Delcev University, Krste Misirkov br. 10-A, Stip, 2000, N. Macedonia} \email{jasmina.buralieva@ugd.edu.mk}

\address{Diana T. Stoeva \newline
Faculty of Mathematics, University of Vienna, 
Oskar-Morgenstern-Platz 1, 1090 Vienna, Austria}
\email{diana.stoeva@univie.ac.at}

\address{Katerina Hadzi-Velkova Saneva \newline
Ss. Cyril and Methodius University in Skopje, Faculty of Electrical Engineering and Information Technologies, Rugjer Boshkovikj 18, 1000 Skopje, N. Macedonia}
\email{saneva@feit.ukim.edu.mk}

\address{Sanja Atanasova \newline
 Ss. Cyril and Methodius University in Skopje, Faculty of Electrical Engineering and Information Technologies, Rugjer Boshkovikj 18, 1000 Skopje, N. Macedonia}
\email{ksanja@feit.ukim.edu.mk}

\keywords{frames, distributions,  generalized asymptotics, Abelian- and Tauberian-type results.\\ \, \\
\indent This is a preprint of the following chapter: J. Veta Buralieva, D. T. Stoeva, K. Hadzi-Velkova Saneva, and S. Atanasova, {\it Asympotic analysis for generalized functions using frames}, accepted for publication in ``Women in Analysis and PDE'' (\url{https://link.springer.com/book/9783031570049}), edited by M. Chatzakou, M. Ruzhansky, and D. Stoeva, part of the sub-series  ``Research Perspectives Ghent Analysis and PDE Center'' of the book series ``Trends in Mahematics'', Birkh\"auser Cham, expected to appear in June 2024. The preprint is posted on arXiv with permission of the Publisher.}
 
\begin{abstract}
This paper is a short overview of the main Abelian- and Tauberian-type results from \cite{SV, ASV, BSSA} regarding the asymptotic analysis of different classes of generalized functions in terms of appropriate frames. The Tauberian-type results provide a comprehensive characterization of the quasiasymptotic and S-asymptotic properties of distributions.
\end{abstract}

\maketitle \numberwithin{equation}{section}
\newtheorem{theorem}{Theorem}[section]
\newtheorem{corollary}[theorem]{Corollary}
\newtheorem{proposition}[theorem]{Proposition}
\newtheorem{lemma}[theorem]{Lemma}
\newtheorem{remark}[theorem]{Remark}
\newtheorem{problem}[theorem]{Problem}
\newtheorem{example}[theorem]{Example}
\newtheorem{definition}[theorem]{Definition}
\allowdisplaybreaks

\section{An introduction to generalized asymptotic analysis}

The theory of \textit{generalized functions} (called also \textit{distributions}) was developed as a result of necessity to find correct mathematical approach for models of various practical processes which were not correctly or not clearly defined by using functions and operations with functions.  There have been many attempts for  generalizing the notion of function, and one of the most significant is the work of  Sobolev in 1936 \cite{Sobolev36}, where he studied the uniqueness of solutions of the Cauchy problem for linear hyperbolic equations.  Also, Bochner in 
  \cite{B1959},  and   Carleman in \cite{C1949}, 
 made an important contribution to the theory of generalized functions by extending the Fourier transform to larger classes of functions. The monograph \cite{S1950} of  Schwartz,  
 which  systematizes the theory of
generalized functions, is based on the theory of linear topological spaces. 
  Among the various approaches to the theory of generalized functions, in the functional approach (the so called Schwartz-Sobolev approach) the theory of distributions is introduced as a part of the duality theory of locally convex spaces. 
Looking at distributions, a main observation is that they do not have point values, contrary to functions which do.
  Point values are a fundamental necessity in most problems
of analysis, so one of the ways to overcome this lack of point value is  through the asymptotic analysis of distributions.

The book of Vladimirov, Drozhzhinov and Zavialov \cite{VDZ} presents the most developed approach to \textit{generalized asymptotics} (i.e., to asymptotic analysis of generalized functions), motivated by theoretical questions in the quantum field theory. In this book the authors developed and applied the theory of quasiasymptotic behavior  in order to analyze asymptotic behavior of  Laplace transforms of tempered distributions. The concept of quasiasymptotic behavior is considered as a natural generalization of the Lojasiewicz concept of point value of a distribution 
\cite{L1}.  Pilipovi\'{c} and his coworkers made an important contribution in this field 
 (see, e.g., \cite{PST, PSV} and references therein). 
  The term S-asymptotic of a distribution was introduced by Pilipovi\'{c} and Stankovi\'{c} \cite{PS1989} (inspired by the shift asymptotics first considered by Schwartz \cite[T. II]{S1950}),
 who furthermore developed a whole theory on S-asymptotic of generalized functions (see, e.g., the book \cite{PSV} and the  
 references therein). 
 Let us note that 
  quasiasymptotics is relevant for the space of tempered distributions, while the S-asymptotics is appropriate for the distributions of exponential type \cite{PST, PSV}.

One way to analyze  the asymptotic behaviors of distributions is through \textit{generalized integral transforms} (i.e., through integral transforms defined on spaces of generalized functions).  
   The book of  Zemanian \cite{Z1968} is the  first systematic monograph on generalized Laplace, Mellin, Hankel, K transform, Weierstrass, the convolution type transform as well as the transforms arising from orthonormal series expansions. Further, generalized Stieltjes,   
  wavelet, short-time Fourier, ridgelet, Radon,
 directional short-time Fourier,  Stockwell transform, and many other generalized integral transforms were defined by different authors. 
Analyzing the asymptotic properties of distributions through the asymptotic behavior of generalized integral transforms, one can obtain two type of results - Abelian and Tauberian. The {Abelian-type result} provides a conclusion for the asymptotic behavior of an integral transform performed on a distribution, under the assumption for some asymptotic behavior of the distribution. 
 The {Tauberian-type result} is the converse of the Abelian-type result, i.e.,  
it provides  a conclusion about the asymptotic behavior of the original distribution, under the assumption for some asymptotic behaviour of the integral transform performed on the distribution. 
 This type of result is more difficult compared to the Abelian one and it usually  requires an  additional condition called Tauberian hypothesis or Tauberian boundedness.  
  Abelian- and Tauberian-type results have already been obtained  for  most of the 
 aforementioned
 generalized integral transforms, we refer e.g. to 
 \cite{BSA, BSpringer,  PST, PTT,PV, PSV, KPSV, SAB, VPR, VDZ}.

In the last several decades, 
 a new research line has been  introduced and developed,  
   in which the asymptotic analysis of a distribution is considered through the asymptotics of the coefficients in its frame expansion. In \cite{Saneva} the author proved several Abelian-type theorems for the asymptotic behaviour of wavelet coefficients of a tempered distribution $f$ assuming the quasiasympotic behaviour of $f$. Motivated by \cite{Saneva}, the authors of \cite{SV} developed a distribution wavelet expansion theory for the space of Lizorkin distributions $\mathcal{S}'_{0}(\mathbb{R})$, and proved several 
 Abelian- and Tauberian-type results
 that relate the quaisasymptotic behavior of distributions to
the asymptotics of wavelet coefficients.
 Further, in \cite{ASV} the authors proved  continuity results for the so-called Gabor  coefficient operator and Gabor synthesis operator on the space of tempered distributions $\mathcal{S}'(\mathbb{R}^{d})$ and the space of distributions of exponential type $\mathcal{K}'_1(\mathbb{R}^{d})$. 
 The characterizations of S-asymptotics of distributions of exponential type is done through Tauberian-type theorems by using Gabor frames. 
  Furthermore, frame expansions through localized frames \cite{GRLF} were extended into 
  the space of tempered distributions in \cite{PS,ps5} 
   and with \cite{SV} and \cite{ASV} motivated consideration of asymptotic analysis of tempered distributions based on localized frames in the preprint \cite{BSSA}. 

In this paper we present some results from \cite{ASV, SV} concerning asymptotic analysis of distributions via Gabor frames and wavelet orthonormal basis, and announce briefly some results from \cite{BSSA} that concern asymptotic analysis of tempered distributions via localized frames.


\section{Preliminaries}

\subsection{Notations} The notation $(f,g)_{\mathcal H}$ stands for the inner product of $f$ and $g$ from a Hilbert space $\mathcal H,$ while $\langle f,\varphi \rangle$ stands for the dual pairing between a distribution $f$ and a test function $\varphi$ so that $(f,\varphi)_{L^2(\R)}=\langle f,\overline{\varphi} \rangle$ when $f$ and $\varphi$ belong to $L^2(\R)$.

\subsection{Spaces of test functions and distributions}\label{owsfd}

The well known \textit{Schwartz space} $\mathcal S(\mathbb R)$ consists of rapidly decreasing smooth functions, that is, $\varphi \in \mathcal S(\mathbb R)$ if  $\varphi\in C^{\infty}(\mathbb R)$ and
\begin{equation}\label{nsw}
\rho_{k}(\varphi)=
\sup_{x\in\R, \alpha \leq k }(1+|x|^2)^{k/2}|\varphi^{(\alpha)}(x)|<\infty, \; \forall k\in\N_0.
\end{equation}
The Fr\'echet space topology of $\mathcal S(\mathbb R)$  is generated by the family of seminorms (\ref{nsw}). Its dual $\mathcal S'(\mathbb R)$ is the space of \textit{tempered distributions} \cite{S1950}. 

The \textit{space of highly time-frequency
localized functions}
$\mathcal{S}_{0}(\mathbb{R})=\{\varphi\in{\s}(\R):\int_{-\infty}^{\infty} x^n\varphi(x) {\rm d}x=0,\forall 
n\in\mathbb{N}_{0}\}$
is a closed subspace of
$\mathcal{S}(\mathbb{R})$ equipped with the relative topology inhered from $\mathcal{S}(\mathbb{R})$. Its dual $\mathcal S'_0(\mathbb R)$ is the space of \textit{Lizorkin distributions} and it is canonically isomorphic to the quotient space of $\mathcal S'(\mathbb R)$ by the space of polynomials \cite{Hol}.

 The \textit{Hasumi-Silva test function space} ${\mathcal K}_1(\mathbb{R} )$  consists of those $\varphi\in C^{\infty}(\mathbb{R})$ for which
$ \sup_{x\in {\mathbb{R} }, \ \alpha\leq k}e^{k|x|}|\varphi^{(\alpha)}(x)| <\infty$, $\forall k\in{\mathbb N}_0.$ The elements of ${\mathcal K}_1(\mathbb{R})$ are called exponentially rapidly decreasing smooth functions.  The dual space  ${\mathcal K}'_1(\mathbb{R})$ consists of all \textit{distributions $f$ of exponential type}, i.e., those of the form $f=\sum_{\alpha\leq l}(e^{s |\:
\cdot\: |}f_\alpha)^{(\alpha)}$, where $f_\alpha\in L^{\infty}(\mathbb{R})$  \cite{hasumi}.


\subsection{Orthonormal wavelet basis} 
A function $\varphi\in
L^2(\R)$  is an \textit{orthonormal wavelet} if the set $\displaystyle \{\varphi_{m,n}:\,m,n\in
\Z\}$ is an orthonormal basis of $L^2(\R)$, where $\displaystyle
\varphi_{m,n}(x)=2^{m/2}\varphi(2^m x-n),\,m,n\in \Z$, \cite{Dob,Walter}. 
So, any
$f\in L^2(\R)$ can be written as
$ f=\sum_{m\in \Z}\sum_{n\in \Z}\left(
f,\varphi_{m,n}\right)_{L^{2}(\mathbb{R})} \varphi_{m,n}$ with convergence in
$L^2(\R)$-norm. This series representation of $f$ is called a \textit{wavelet expansion}. 
 The \textit{wavelet coefficients} of $f$ with respect to the orthonormal wavelet
$\varphi$  will be denoted by $c_{m,n}^{\varphi}(f)$,
\begin{equation*} c_{m,n}^{\varphi}(f)=\left(
f,\varphi_{m,n}\right)_{L^{2}(\mathbb{R})}=\int_{-\infty}^{\infty}f(x){\bar\varphi}_{m,n}(x){\rm
d}x,\,\,m, n\in \Z\: .\end{equation*} 

It is well known that every orthonormal
wavelet from ${\s}(\R)$ must belong to the space ${\s}_0(\R)$
\cite[Cor.3.7]{HW}.


\subsection{Frames} 

The frame concept (introduced in \cite{DS} for Hilbert spaces) extends the concept of an orthonormal basis, gives much more freedom for construction of frame elements in comparison to the orthonormal basis restrictions (e.g. redundancy is possible),  
and still guarantees perfect and stable reconstruction of all the elements of the space, which makes the concept very useful for applications. 
 A sequence $(u_n)_{n\in I}$  with elements from a Hilbert space $\h$  ($I$ being a countable index set)
 is called a \textit{frame for $\h$} if  there exist positive constants $A$ and $B$ so   that $A\|h\|^2 \leq \sum_{n\in I} |{(h, u_n)_{\h}} | ^2 \leq B\|h\|^2$ for every $h\in\h$. Given a frame $(u_n)_{n\in I}$ for $\h$, there always exists a frame $(v_n)_{n\in I}$ for $\h$ such that every $h\in\h$ can be reconstructed from the \textit{frame coefficients} ${(h,u_n)_{\h}}, n\in I,$ 
using the frame $(v_n)_{n\in I}$ (called a \textit{dual frame} of $(u_n)_{n\in I}$) as follows: $h=\sum_{n\in I} {(h, u_n)_{\h}} v_n$. This reconstruction formula is also refered as \textit{frame expansion} of the elements of $\h$. A frame for $\h$ that is at the same time a Schauder basis for $\h$ is called a Riesz basis for $\h$. 
The unique biorthogonal sequence to a Riesz basis $(r_n)_{n=1}^\infty$ for $\h$  is also a Riesz basis for $\h$ and it will be denoted by $(\widetilde{r}_n)_{n\in I}$. 

Frames of specific structure find application in many fields. In particular, Gabor  frames play significant role in signal and image processing. 
Given $\varphi\in L^2(\mathbb R)$ and positive constants $\alpha$ and $\beta$, a \textit{Gabor frame for} {$L^2(\mathbb R)$} is a frame for  $L^2(\mathbb R)$ of the form 
$G(\varphi,\alpha,\beta)=\{e^{2\pi {\rm i}m\beta \cdot} \varphi(\cdot - n\alpha)\}_{m,n\in\mathbb Z}$.
  For more on general frame theory, as well as on Gabor frames, we refer e.g. to the 
books \cite{Olebook,FS1}.

Several types of localization of frames were introduced in the literature by {G}r{\"o}chenig 
and by 
Balan, Casazza, Heil, and Landau. 
We use the notion of polynomially localized frames introduced in \cite{GRLF}.  
A frame $(u_n)_{n=1}^\infty$ for $\h$ is called \textit{polynomially localized  with decay $\gamma > 0$} (in short, \textit{$\gamma$-localized}) \textit{with respect to a Riesz basis $(r_n)_{n=1}^\infty$ for $\h$} if 
there exists a constant $C_\gamma$ such that
$$ \max \{|{( u_m, r_n)_{\h}|, |(u_m, \widetilde{r}_n)_{\h}}|\}
\leq C_\gamma(1+|m-n|)^{-\gamma}, \ \forall m,n\in \mathbb N. $$ 
An important property of such frames is that the frame coefficients of a function lie not just in $\ell^2$, but also in a class of associated sequence spaces (which is helpful for determining certain properties of the function like decay, smoothness, and other)  
and the frame expansion on the Hilbert space extends on a class of associated Banach spaces \cite{GRLF} and on the space of tempered distributions \cite{ps5}.


\section{Asymptotic Behavior of Distributions}\label{owabd}

In this section we briefly explain the asymptotic notions for
 distributions that will be used in the next sections.

\subsection{Quasiasymptotic behavior of distributions} 
A measurable real-valued function $L$ 
defined and positive on an interval  $(0,A]$ (resp. $ [A,
\infty)), $ where $A$ is a positive constant, 
is called \textit{slowly varying at the origin
(resp. at infinity)} if
$$
\lim_{\varepsilon\rightarrow0^+}\frac{L(a\varepsilon)}{L(\varepsilon)}=1
\,\,\,\Big(\mbox{resp. }
\lim_{\lambda\rightarrow\infty}\frac{L(a\lambda)}{L(\lambda)}=1\Big)\
\ \text{for each}\ a>0\ .
$$
We refer to \cite{Seneta} for  properties of the slowly varying functions.

Let $L$ be a slowly varying function at the origin (resp. at infinity). It is said that
the distribution $f\in{\s}'(\R)$ has \textit{quasiasymptotic behavior of
degree $\alpha\in \R$ at the point $x_0\in\mathbb{R}$  (resp. at infinity)} \textit{with respect
to $L$} if there exists $g\in\mathcal{S}'(\R)$ such that for each
$\varphi\in\mathcal{S}(\R)$
\begin{equation} \label{equation14} \lim_{\varepsilon \rightarrow 0^+}\Big\langle
\frac{f(x_0+\varepsilon \cdot)}{\varepsilon^\alpha
L(\varepsilon)},\,\varphi\Big\rangle=\langle g,\varphi\rangle  \,\,\Big (\mbox{resp.} \, \lim_{\lambda \rightarrow \infty} 
\Big\langle
\frac{f(\lambda\cdot)}{\lambda^\alpha L(\lambda)},\,\varphi\Big\rangle\Big ).
\end{equation}

 We use the following convenient notation  
$ f(x_0+\varepsilon x)\sim\varepsilon^\alpha L(\varepsilon)
g(x)\,\,\mbox{ as}\,\,\varepsilon\rightarrow 0^+\,\,$ $  (\mbox{resp.} \, \,  f(\lambda x)\sim\lambda^\alpha L(\lambda) g(x) \,\,\mbox{ as}\,\,\lambda \rightarrow \infty\,) \, \,{\mbox
{in}}\,\, {\s}'(\R) $
 with the meaning of (\ref{equation14}). 
Under the above setting, 
 one can prove that $g$ must be homogeneous with degree of homogeneity $\alpha$, namely, $g(ax)=a^\alpha g(x)$, for each $a>0$ \cite{PSV}. 
 For more details on quasiasymptotic theory, we refer to
\cite{PST,VDZ, PSV}.

We also consider quasiasymptotics in $\mathcal{S}'_{0}(\R)$. 
The relation $f(x_0+\varepsilon x)\sim\varepsilon^\alpha
L(\varepsilon) g(x)$ as $\varepsilon\rightarrow
0^+$ in $\mathcal {S}'_{0}(\R)$ means that 
the respective limit in 
\eqref{equation14} is satisfied just for each
$\varphi\in\mathcal{S}_{0}(\R)$; and analogously for quasiasymptotics
at infinity in $\mathcal{S}'_{0}(\R)$.

\subsection{$S$-asymptotic behavior of distributions}
The natural framework for $S$-asymptotic behavior of distributions is the space of distributions of exponential type $\mathcal{K}'_{1}(\mathbb{R})$ \cite[Chap.~1]{PSV}. Let   $c:\mathbb{R}\to (0,\infty)$ be a measurable function. 
It is said that $f\in\mathcal{K}'_1(\mathbb{R})$ has \textit{$S$-asymptotic behavior} (\textit{$S$-asymptotics}) \textit{with respect to $c$} if there is $g\in\mathcal{K}'_1(\mathbb{R})$ such that  for each
$\varphi\in\mathcal{K}_1(\R)$ 
\begin{equation}
\label{Seq1}\lim_{h\to\infty} \Big\langle
\frac{f(\cdot+h)}{
c(h)},\,\varphi\Big\rangle=\langle g,\varphi\rangle  .
\end{equation}

 As shown in \cite{PSV}, relation (\ref{Seq1}) forces the  distribution $g$ and the function $c$ to have the form $g(x)=C e^{b x}$ and $c(h)=e^{bh}L(e^h)$, for some $C\in\mathbb{R}$, $b\in \mathbb{R}$, and  
slowly varying function $L$ at infinity.

\section{Asymptotic results for {distributions of exponential type} through Gabor frames }

Let $G(\varphi,\alpha,\beta)=\{e^{2\pi {\rm i}m\beta \cdot} \varphi(\cdot - n\alpha)\}_{m,n\in\mathbb Z}$  be a Gabor frame with window $\varphi\in \mathcal{K}_{1}(\mathbb{R})$. The \textit{Gabor frame coefficients} of $f\in\mathcal{K}_{1}'(\mathbb{R})$ are 
$$\langle f, e^{-2\pi {\rm i}m\beta \cdot}\, \overline{\varphi(\cdot - n\alpha)}\rangle =V_{\varphi}f(\alpha n, \beta m),$$
where 
 $V_{\varphi} f(x,\xi )=\int _{-\infty}^\infty
f(t)\overline{\varphi(t-x)}e^{-2\pi i\xi  t} dt, x,\xi \in {\mathbb{R}},$
is the {short-time Fourier transform} (STFT)  of a function $f$ with respect to a window function $\varphi$.

The next theorems can be interpreted as Tauberian-type results   
for the $S-$ asymptotics of distributions through the growth estimates of Gabor frame coefficients. 

As noticed at the end of Section 3.2,  
in the considered space $\mathcal{K}_{1}'(\mathbb{R})$, 
expressions for the functions $g$ and $c$ in the S-asymptotics (\ref{Seq1}) are explicitely determined. Rewriting \cite[Theorem 2]{ASV} for this case, the following holds:

\begin{theorem} \cite{ASV} \label{ttGF1} Let $G(\varphi,\alpha,\beta)$ be a Gabor frame with window $\varphi\in \mathcal{K}_{1}(\mathbb{R})$.  
Then, $f\in\mathcal{K}_{1}'(\mathbb{R})$ has the $S$-asymptotic behavior 
 $(\ref{Seq1})$ 
if and only if there are slowly varying function $L$ and $b\in\R$ so that   
\begin{equation*}\label{limit2} \lim_{x\to \infty }e^{2\pi i \beta m  x}\frac{V_{\varphi}f(x,\beta m)}{e^{bx}L(e^x) }=:a_m\in\mathbb{C}\quad \mbox{exists for every } m\in\mathbb{Z}\end{equation*} 
and 
 there is $\tau\in\mathbb{R}$ such that
\begin{equation*} 
\sup_{(x,m)\in\mathbb {R}\times \mathbb{Z}}\frac{|V_\varphi f(x,\beta m)|}{e^{bx}L(e^x)(1+|m|)^{\tau}}<\infty.
\end{equation*} 

Under the above equivalence, the limit function $g$ in the S-asymptotic relation $(\ref{Seq1})$ has the form $g(x)=C e^{bx},$ 
where $C$ and $b$ are related to $a_m$ as follows: 
$a_m=C \overline{\widehat{\varphi}\left(-\beta m+ib/(2\pi)\right)}$, $m\in\mathbb{Z}.$

\end{theorem}

For the specific case of non-decreasing functions one has the following result concerning their S-asymptotics:

\begin{theorem}\cite{ASV} \label{ttGF2} Let $f$ be a positive non-decreasing function on $[0,\infty)$, $let$ $L$ be a slowly varying function, and
let $G(\varphi,\alpha,\beta)$ be a Gabor frame with nonnegative window $\varphi\in \mathcal{K}_{1}(\mathbb{R})$.
 Suppose that the limits
\begin{equation*}
\label{eq01}\lim_{x\to\infty}\frac{e^{2\pi i m x}}{e^{b x}L(e^{x})} \int _{0}^{\infty} f(t)\varphi(t-x)e^{-2\pi i mt} \ dt=: a_m
\end{equation*} exist for all $m\in\mathbb{Z}$. Then,
\begin{equation*}
\lim_{x\to\infty} \frac{f(x)}{e^{b x}L(e^{x})}= \frac{a_0}{\int_{-\infty}^{\infty}\varphi(t)e^{b t} dt}\ .
\end{equation*}
\end{theorem}


\section{Asymptotic results for tempered distributions through wavelets}

In this section we present some Tauberian-type results from \cite{SV} that give complete 
characterization of quasiasymptotics  in terms of wavelet coefficients.

Let $\varphi\in
{\mathcal{S}}_0(\mathbb{R})$ be an orthonormal wavelet and let
$f\in\mathcal{S}'(\mathbb{R})$. The wavelet coefficients of
$f(x_{0}+\varepsilon\:\cdot)$, $x_0\in\mathbb{R}$, are denoted by 
\begin{equation*} 
c_{m,n}^{\varphi}(f;\varepsilon,x_0)=\langle f(x_0+\varepsilon
\cdot),\bar{\varphi}_{m,n}\rangle\,.
\end{equation*}  

The next theorem provides a necessary and sufficient condition for the existence of  quasiasymptotics in the space $\mathcal{S}'_{0}(\mathbb{R})$. 
\begin{theorem} \cite{SV}\label{owtheorem} Let $\varphi\in
{\mathcal{S}}_0(\mathbb{R})$ be an orthonormal wavelet and let $L$ be slowly varying at the origin. For $f\in {\mathcal S}'_{0}(\mathbb{R})$,
the following two conditions,
\begin{enumerate}
\item [(i)] the limit  $ \lim_{\varepsilon\rightarrow
0^+}\frac{c_{m,n}^{\varphi}(f;\varepsilon ,
x_0)}{\varepsilon^{\alpha}
L(\varepsilon)}<\infty
$ exists for each
$m,n\in \mathbb{Z}$, and 
\item[(ii)] there exist positive constants $\beta,\gamma$, and $C$ such that 
$$
|
\,c_{m,n}^{\varphi}(f;\varepsilon ,x_0)| \leq
C\varepsilon^{\alpha} L(\varepsilon)(1+|n|)^{\beta}\left(2^m+\frac{1}{2^m}\right)^{\gamma}\ , 
$$
 for
all $m,n\in \mathbb{Z}$ and $0<\varepsilon\leq 1$,
\end{enumerate}
are necessary and sufficient 
 for the existence of a distribution $g\in \mathcal{S}'_0(\mathbb{R})$ 
 such that 
$$f(x_0+\varepsilon
x)\sim\varepsilon^{\alpha}L(\varepsilon)g(x)\,\,\,{\mbox
\,{as}}\,\,\varepsilon\rightarrow 0^+\,\,\,{\mbox
in}\,\,\,{\mathcal{S}}'_0(\mathbb{R})\,.$$
\end{theorem}

The next statements are Abelian- and Tauberian-type results, relating the quasiasymptotic behavior of tempered distributions 
with the asymptotics of wavelet coefficients. 

\begin{proposition} \cite{SV} 
\label{owp3} Let $\varphi\in {\s}_0(\R)$ be an orthonormal wavelet. 
If  $f\in {\s}'(\R)$  has the quasiasymptotic behavior \eqref{equation14} 
at the point $x_0$, 
then for
each $m,n\in\Z$,
$
\lim_{\varepsilon\rightarrow
0^+}\frac{c_{m,n}^{\varphi}(f;\varepsilon,x_0)}{\varepsilon^{\alpha}L(\varepsilon)}=c_{m,n}^{\varphi}(g).$ 
\end{proposition}

\begin{theorem} \cite{SV}\label{theorem5}
Let $\varphi\in {\s}_0(\R)$ be an orthonormal wavelet. Let $f\in
{\s}'(\R)$ and 
$\alpha >0,\,\alpha\notin {\N}$. The conditions {\rm (i)} and {\rm (ii)} from Theorem \ref{owtheorem}
are necessary and
sufficient for the existence of a polynomial $p$ of degree less
than $\alpha$ and a homogeneous distribution $g$ of degree
$\alpha$ such that
\begin{equation*}\label{4.3'} f(x_0+\varepsilon x)=p(\varepsilon
x)+\varepsilon^{\alpha}L(\varepsilon)g(x)+o(\varepsilon^{\alpha}L(\varepsilon))\,\,{\mbox
{as}}\,\,\,\varepsilon\rightarrow 0^+\,\,\,{\mbox
in}\,\,\,{\s}'(\R)\ .\end{equation*}
\end{theorem}

Similar assertions to the previous ones hold for quasiasymptotics at infinity.

\section{Asymptotic results for tempered distributions through localized frames}

In this section we announce some results from \cite{BSSA} that concern Abelian- and Tauberian-type results for quasiasymptotic behaviour of tempered distributions through  polinomially localized frames with respect to the Hermite basis. 
Recall, the Hermite orthonormal basis
$(h_n)_{n=1}^\infty$ of $L^2(\mathbb R)$ is defined as $ h_n = {\rm h}_{n-1}$, $n\in\mathbb N$, where
$${\rm h}_n(t) = (2^{n} n!\, \sqrt{\pi})^{-1/2}  \, (-1)^{n} e^{t^2/2} \frac{d^{n}}{dt^{n}} (e^{-t^2}),\,n\in\mathbb N_0.$$

\begin{theorem} \label{owtheorem1} 
 Assume that $(e_n)_{n=1}^\infty$  is a sequence with elements from $\mathcal{S}(\mathbb R)$ which is a frame for $L^2(\mathbb R)$
and which is  polynomially localized with respect to the Hermite basis $(h_n)_{n=1}^\infty$ with decay $s$ for every $s\in\mathbb N$.
Let $L$ be slowly varying at the origin, $\alpha \in \mathbb R$, and $x_0\in\mathbb R$.  
 For $f\in {\mathcal S}'(\mathbb R)$, there exists  a tempered distribution $g$ such that 
(\ref{equation14})  holds   at $x_0$ if and only if the following two conditions hold: \begin{enumerate} \item [(i)] the limit 
$\lim_{\varepsilon\rightarrow 0^+}\langle \frac{ f(\varepsilon x + x_0)}{\varepsilon^{\alpha} L(\varepsilon)}, e_n \rangle$ exists for each $n\in \mathbb N$;
 
 \item[(ii)] there exist $k\in\mathbb N_{0}$ and  $\varepsilon_0 \in (0,\infty)$ such that  
$$ \sup_{\varepsilon\in(0,\varepsilon_0)} \sup  _{n\in\mathbb N} \frac{| \, \langle f (\varepsilon x+ x_0),  e_n \rangle |} { \varepsilon^{\alpha} L(\varepsilon)n^k} <
 \infty.$$

 \end{enumerate}
 \end{theorem}

\begin{theorem} \label{owtheorem2} 
Assume that the conditions of 
Theorem \ref{owtheorem1} are satisfied.
Let $L$ be slowly varying at infinity and  $\alpha \in \mathbb R$.  
 For $f\in {\mathcal S}'(\mathbb R)$, there exists  a tempered distribution $g$ such that 
 (\ref{equation14}) holds at infinity if and only if the following two conditions hold: \begin{enumerate} \item [(i)] the limit   
$\lim_{\lambda\rightarrow \infty}\langle \frac{ f(\lambda x)}{\lambda^{\alpha} L(\lambda)},  e_n \rangle$
exists for each $n\in \mathbb N$; 
 \item[(ii)] there exist $k\in\mathbb N_{0}$ and $\lambda_0>0$ such that  
$$  \sup_{\lambda\geq \lambda_0}  \sup  _{n\in\mathbb N} \frac{| \, \langle f(\lambda x),  e_n \rangle |} { \lambda^{\alpha} L(\lambda)n^k} <
 \infty. $$

 \end{enumerate}
\end{theorem}

\section*{Acknowledgements}

The second author is supported from the Austrian Science Fund (FWF) through  Project P 35846-N ``Challenges in Frame Multiplier Theory''.

\end{document}